\documentclass{amsart}

\usepackage{a4,xy,amssymb}

\xyoption{all}

\def\cM{{\mathcal M}}

\def\C{{\Bbb C}}
\def\Rl{{\Bbb R}}

\def\Z{{\Bbb Z}}

\def\eop{\hspace*{\fill}$\square$}
\def\implies{\Rightarrow}

\def\to{\rightarrow}

\def\<{\langle}
\def\>{\rangle}

\theoremstyle{definition}
\newtheorem{defn}{Definition}[section]
\newtheorem{thm}{Theorem}[section]

\newtheorem{lemma}[thm]{Lemma}
\theoremstyle{remark}

\newtheorem*{prf}{Proof}
\newcounter{kd}

\newcommand{\grp}[1]{\ensuremath{\mathsf{#1}}}

\newcommand{\iss}{\ensuremath{\mathsf{iss}}}
\newcommand{\GL}{\ensuremath{\mathsf{GL}}}
\newcommand{\Repr}{\ensuremath{\mathsf{Rep}}}
\newcommand{\Proj}{\ensuremath{\mathsf{Proj~}}}
\newcommand{\Mat}{\ensuremath{\mathsf{Mat}}}
\newcommand{\Hom}{\ensuremath{\mathsf{Hom}}}

\newcommand{\se}[1]{\begin{equation*}\begin{split}#1\end{split}\end{equation*}}

\title{Smooth Character Varieties for Torus Knot Groups}
\newcommand {\dss }{/\!/}

\newcommand {\vertex}[1]{*+[o][F-]{#1}}
\newcommand {\Mss}{\mathsf{M}^{ss}}
\newcommand {\PP}{\mathbb P}
\newcommand {\Qpq}{Q_{pq}}
\newcommand{\vtx}[1]{*+[o][F-]{\scriptscriptstyle #1}}

\author{Jan Adriaenssens}
\address{Jan Adrianssens\\ Universiteit Antwerpen (RUCA) \\ B-2610 Antwerp (Belgium)}
\email{jan.adriaenssens@ua.ac.be}

\author{Raf Bocklandt}
\address{Raf Bocklandt\\ Universiteit Antwerpen (UIA) \\ B-2610 Antwerp (Belgium)}
\email{raf.bocklandt@ua.ac.be}

\author{Geert Van de Weyer}
\address{Geert Van de Weyer\\ Universiteit Antwerpen (UFSIA) \\ B-2610 Antwerp (Belgium)}
\email{geert.vandeweyer@ua.ac.be}

\begin{document}
\bibliographystyle{alpha}

\begin{abstract}
Semisimple representations of the free product $\Z_p*\Z_q$ determine $\theta$-semistable representations of a specific quiver $Q_pq$. The dimension vectors of $\theta$-stable representations of this quiver were classified in \cite{AdLeB}. In this paper we classify the moduli spaces $\Mss_\alpha(\Qpq,\theta)$ which are smooth projective varieties.
\end{abstract}
\maketitle

\section{Introduction}
Consider a cilinder with $q$ line segments on its surface, equidistant and parallel to its axis. If the ends of this cilinder are identified with a twist $2\pi\frac{p}{q}$ where $p$ is an integer relatively prime to $q$, one obtains a single curve on the surface of a torus. Such a curve is called a torus knot, and is denoted by $K_{p,q}$. The fundamental group of the complement $\Rl^3\backslash K_{p,q}$ is called the $(p,q)$-torus knot group and is equivalent to the group
$$\langle x,y \mid x^p = y^q\rangle.$$
The center of this torus knot group is generated by the element $y^q$, so the quotient of a torus knot group with its center is equivalent to the free product $\Z_p*\Z_q$. If one wants to study irreducible representations of such a torus knot group, it suffices to study the representation theory of the quotient, and hence of $\Z_p*\Z_q$.
In \cite{AdLeB}, Adriaenssens and Le~Bruyn show that one can reduce the complex representation theory of the free product of two finite cyclic groups to the representation theory of a certain bipartite quiver.

The equivalence between representations of $\Z_p * \Z_q$ and representations of quivers is achieved in the following way. Consider a complex representation $V$ of the free product. By looking only at the action of $\Z_p$, one can decompose the 
vectorspace $V$ into a direct sum of eigenspaces $V_{\xi^1} \oplus \dots V_{\xi^p}$ where $\xi$ is a $p^{th}$ root of unity. Doing the same regarding to $\Z_q$, we obtain a double decomposition:
$$\xymatrix{
V_{\xi^1} \oplus \dots \oplus V_{\xi^p} \ar[rr]^{\cong} && V \ar[rr]^{\cong} && W_{\eta^1} \oplus \dots \oplus W_{\eta^q} \ar@{->>}[d]\\
V_{\xi^i} \ar@{^(->}[u]\ar[rrrr]_{m_{ij}}&& && W_{\eta^j}
}
$$

So the canonical situation is that we have $p+q$ vectorspaces and a linear map from each of the $p$ first spaces to each of the $q$ last. This is in fact a representation of the folowing quiver:
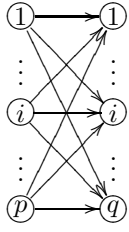
$$
\xymatrix
{
\vertex{1} \ar@{}[d]|{\vdots}\ar[r] \ar[rd] \ar[rdd] & \vertex{1} \ar@{}[d]|{\vdots}\\
\vertex{i} \ar@{}[d]|{\vdots}\ar[ur]\ar[r]\ar[dr]    & \vertex{i} \ar@{}[d]|{\vdots}\\
\vertex{p} \ar[ruu] \ar[ru] \ar[r] & \vertex{q}
}
$$

The only restriction on the maps is that they must add up to an invertible map $M$ between $V_{\xi^1} \oplus \dots \oplus V_{\xi^p}$ and $W_{\eta^1} \oplus \dots \oplus W_{\eta^q}$, because all the maps are in fact restrictions of the indentity on $V$. This condition is necessary and sufficient. If we define the dimension vector of a $\Z_p*\Z_q$-representation as the vector
$$
\alpha := (\grp{Dim} V_{\xi^1}, \dots, \grp{Dim} V_{\xi^p};\grp{Dim} W_{\eta^1}, \dots, \grp{Dim} W_{\eta^q}),
$$

we can say that there is an equivalence of categories between the category \mbox{$\grp{Rep_\alpha}\Z_p*\Z_q$}, containing the representations with dimension vector $\alpha$ and the Zariski open subset $U$ of $\Repr_\alpha Q$, consisting of the $\alpha$-dimensional representations of the quiver for which the block matrix $M$ is invertible. The action of $\grp{GL_n}$ on $\grp{Rep_\alpha}\Z_p*\Z_q$ translates itself into an action of $\grp{GL_\alpha} = \prod_i\grp{GL}_{\alpha_i}$ on $\Repr_\alpha Q$. So classifying the representation classes in $\grp{Rep_\alpha}\Z_p*\Z_q$ is the same thing as classifying the orbit of the $\grp{GL_\alpha}$-action in $U$. Doing this we will obtain that $\grp{iss}_\alpha \Z_p*\Z_q \equiv U\dss \grp{GL_\alpha}$ is an affine variety containing the semisimple representation classes of $\Z_p*\Z_q$. 

A geometrically more appealing approach to study this affine variety is to look at a certain projective closure of this variety, the moduli space of $\alpha$-dimensional $\theta$-semistable representations of the quiver.
\begin{defn}
Let $\theta$ be the following vector $(-1,\dots,-1;1,\dots, 1) \in \Z^{p+q}$. An $\alpha$-dimensional representation of the quiver $Q$ is said to be \emph{$\theta$-semistable} if and only if
\begin{itemize}
\item $\theta \cdot \alpha =0$
\item For every subrepresentation with dimension vector $\beta$, $\theta \cdot \beta \ge 0$
\end{itemize}
If we can put a strict inequality in the last item, the representation is called \emph{$\theta$-stable}.
\end{defn}

One can easily check that every representation in $U$ is in fact a $\theta$-semistable. Indeed, if there is a subrepresentation with $\theta \cdot \beta<0$, the big map $M$ maps a subspace $V_1' \oplus \dots \oplus V_p'$ onto a subspace $W_1' \oplus \dots \oplus  W_q'$ of smaller dimension so $M$ is definitely non-invertible. This implies that $U$ is a open and dense subset of the $\theta$-semistable representations of $Q$.

A second property of the $\theta$-semistable representations is that a closed $\grp{GL_\alpha}$-orbit in $U$ is also closed in the variety of $\theta$-semistable representations. If this would not be the case there would be $(X;Y) \in \grp {Mat}_{\alpha}(\C)$ such that $XMY$ is not invertible and $M$ is. This implies that either $X$ has a kernel $V_1' \oplus \dots \oplus V_p'$ or $Y$ has an image $W_1' \oplus \dots \oplus  W_q'$. In the first case there is a subrepresentation of $XMY$ with dimension vector $(\grp{Dim}V_1',\dots,\grp{Dim}V_p; 0,\dots ,0)$ in the second case there is one with dimension vector $(\grp{Dim}V_{\xi^1},\dots,\grp{Dim}{V_\xi^p};\grp{Dim} W_1',\dots ,\grp{Dim} W_q')$. Both will give a negative number when multiplied by $\theta$ so if $XMY$ is not invertible, it is also not $\theta$-semistable.

The two properties mentioned above, enable us to view the quotient variety $U\dss \grp{GL_\alpha}$ as an open dense subvariety of $\Mss_\alpha(Q,\theta) := \Repr^{ss}_\alpha(Q,\theta)\dss \grp{GL_\alpha}$ and we have thus the following
diagram.
$$\xymatrix{
U \ar@{^(->}[r] \ar@{->>}[d] &\Repr^{ss}_\alpha (Q,\theta) \ar@{->>}[d]\\
U\dss \GL_\alpha \ar@{^(->}[r] &M^{ss}_\alpha(Q,\theta)
}
$$
This diagram indicates that to study the representations of $\Z_p*\Z_q$ one could first try to study the moduli space
$M^{ss}_\alpha(Q,\theta)$.

From now on we're going to work only with $\theta$-semistable representations of the quiver $Q$. So it's time to fix some notation. The vector spaces on each vertex will be denoted by $V_i, p = 1,\dots, p$ for the left vertices of the quiver and $W_i, p = 1,\dots, q$ for the right ones.

The semistability implies that the dimensionvector is of the following form
$$
\alpha := (a_1,\dots,a_p; b_1,\dots,b_q)\text{ where }\sum_{i=1}^p a_i = \sum_{i=1}^q b_i:=n
$$
If we look at the Euler form of the quiver $Q$, i.e. the matrix with entries
$$
[\chi_Q]_{ij} := \delta_{ij} - \#\{\text{arrows from vertex $i$ to $j$ }\}.
$$
We can decompose it to a block matrix of the following form
$$
\left[
\begin{matrix}
1 &       & 0 & -1   &\cdots & -1\\
  &\ddots &   &\vdots&       &\vdots\\
0 &       & 1 &-1    &\cdots &-1\\
  &       &   & 1    &       & 0\\
  &  0    &   &      &\ddots &  \\
  &       &   & 0    &       & 1\\
\end{matrix}
\right]
$$
Now consider two dimension vectors $\alpha_1$ and $\alpha_2$. One can easily compute their image under the Euler form:
$$
\chi_Q(\alpha_1,\alpha_2) = \alpha_1\cdot \alpha_2-n_1n_2,
$$
where $n_1$ resp.~$n_2$ equals the sum of the first $p$ resp.~the last $q$ entries of the dimension vector.
For the remainder of the paper, $n$ shall always equal the sum of the entries of a dimension vector considered.

The last convention we make is that we shall write elements of $\GL_\alpha$ in  this way:
$$
g := (g_1,\dots,g_p;h_1,\dots,h_q)~g_i \in \GL_{a_i},h_i \in \GL_{b_i}.
$$

\section{The local structure of the moduli space $\Mss_\alpha(Q,\theta)$}
In \cite{King94}, King showed that $\Mss_\alpha(Q,\theta)$ has the structure of a projective variety. Algebraically it corresponds to the graded ring of semi-invariant functions with character 
$$
\chi_\theta: \grp{GL_\alpha} \to \C^*: g \mapsto (\det g_1 \cdot \dots \cdot \det g_p)^{-1}(\det h_1 \cdot \dots \cdot \det h_p).
$$

If we extend the space $\Repr_\alpha Q$ to $\Repr_\alpha Q\oplus \C$ together with an extended action
$$
\forall (V,c) \in \Repr_\alpha Q\oplus \C  (V,c)^g = (V^g,c\chi_\theta(g)^{-1}),
$$
The ring of polynomials over $\Repr_\alpha Q\oplus \C$ is of the form $\C[\Repr_\alpha Q][t]$ and becomes a graded ring by defining
$$
\deg t =1,\forall f \in\C[\Repr_\alpha Q]:  \deg f =0 
$$
We can consider the subring of invariant polynomial functions on $\Repr_\alpha Q\oplus \C$, which is also graded in the same way
$$
\C[\Repr_\alpha Q\oplus \C]^{\GL_\alpha}:= \{\sum_i f_it^i|\forall g \in \GL_\alpha :f_i\circ g=\chi_\theta^{i}(g)f_i\} 
$$
This graded ring corresponds to a projective variety, $\Proj \C[\Repr_\alpha Q\oplus \C]^{\GL_\alpha}$, consisting of the graded-maximal ideals not containing the positive part $\C[\Repr_\alpha Q\oplus \C]^{\GL_\alpha}_{\deg>0}$. If $M$ is a graded-maximal ideal in $\C[\Repr_\alpha Q\oplus \C]^{\GL_\alpha}$ then it is contained in a maximal ideal of the ring 
$\C[\Repr_\alpha Q\oplus \C]$ which corresponds to a couple $(V,c)$. Moreover if $M$ doesn't contain the positive part $c$ is definitely not zero and there exists at least one $ft^n \in \C[\Repr_\alpha Q\oplus \C]^{\GL_\alpha}$ such that $f(V)\ne 0$, such an $f$ is called a \emph{semi-invariant of $V$}. Vice versa if $V$ is a representation such that there exists a
semi-invariant then 
$$
M_{V}:=\{\sum_i f_it^i \in \C[\Repr_\alpha Q\oplus \C]^{\GL_\alpha}| f_i(V) = 0\},
$$
will be a maximal-graded ideal not containing the positive part.

A method for constructing these semi-invariants was discovered independently by Schofield and Van~den~Bergh in \cite{SVdB}, Derksen and Weyman in \cite{DW} and Domokos and Zubkov in \cite{DZ}. Take two diagonal matrices $A \in \Mat_{p\times p}\C Q$ and $B \in \Mat_{q\times q}\C Q$ such that the diagonal elements are the vertices of the quiver. Consider a matrix $\cM \in A^{\oplus n} \Mat_{np \times nq}\C Q B^{\oplus n}$, the entry $\cM_{ij}$ is now a linear combination of paths from $v_{i \mod p}$ to $w_{j \mod q}$.
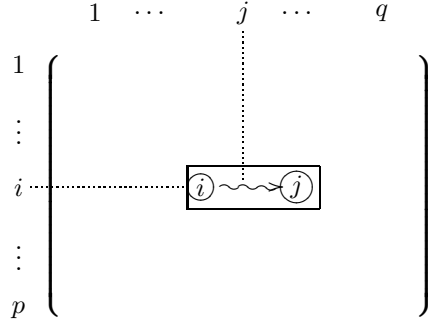
$$
\xymatrix@R-1.5pc@C-1.5pc{
       &       & 1 & \dots &  & j\ar@{.}[ddd] & \dots && q &\\
1      &       &   &       &                &&       &&   &\\
\vdots &       &   &       &                &&       &&   &\\
i\ar@{.}[rrrr]&&   &       &  \vertex{i}\ar@{~>}[rr]&&\vertex{j}&    &   &\\
\vdots &       &   &       &                &&       &&   &\\
p      &       &   &       &                &&       &&   &
\save "2,2"."6,2"*+<-.1pc>\frm{(}
\restore
\save "2,10"."6,10"*+<-.1pc>\frm{)}
\restore
\save "4,5"."4,7"*+<+.3pc>[F-]\frm{}
\restore
}
$$
Using a representation $V$ of $Q$ we can map each $\cM_{ij}$ to a linear map in $\Hom_\C (V_{i \mod p},W_{i \mod q})$. Putting all those maps together we get the linear map 
$$
\cM_V : (\bigoplus_{i=1}^{p} V_{i})^{\oplus n} \to (\bigoplus_{i=1}^{n_2} W_i)^{\oplus n}
$$
If the dimensions of the source and target of these map are the same, we can take the determinant of this map. This determinant varies under de action of $\GL_\alpha$ as:
\se{
\det(\cM_{V^g}) &= \det \left( (\bigoplus_{i=1}^{p} g_i^{-1})^{\oplus n}\cM_V   (\bigoplus_{i=1}^{q} h_i)^{\oplus n} \right)\\
&=\prod_{i=1}^p \det g_i ^{-n}\prod_{i=1}^q \det h_i ^{n} \det(\cM_V)\\
&=\chi_\theta(g)^n \det(\cM_V)\\
}
So $f_{\cM} : V \mapsto \det \cM_V$ is a semi-invariant of order $n$. We now have a way to construct semi-invariants and one can even prove that those semi-invariants generate all invariants. This observation leads to the following lemma
\begin{lemma}
If $V$ is a $\theta$-semistable then there is a matrix-semi-invariant $\cM$ such that
$$
f_\cM(V) \ne 0
$$
\end{lemma}

In order to determine which moduli spaces are smooth projective varieties, we will use a result by Le~Bruyn and Procesi \cite{ProLeB} determining the local structure around a point $V\in \Mss_\alpha(Q,\theta)$. Let
$$V = S_1^{\oplus a_1}\oplus\dots\oplus S_k^{\oplus a_k}$$
where $S_i$ is a $\theta$-stable representation of dimension vector $\alpha_i$. The \emph{local quiver} $Q_V$ of this representation is a quiver on $k$ vertices (corresponding to the $k$ distinct terms in the decomposition) with the number of arrows between vertices $i$ and $j$ determined by
$$\delta_{ij} - \chi_Q(\alpha_i,\alpha_j).$$
Note that in the case of a bipartite quiver this number equals $\delta_{ij} + n_in_j - \alpha_i.\alpha_j$. The multiplicities of each term in $V$ yield a dimension vector for this quiver:
$$\beta = (a_1,\dots, a_k).$$
This local quiver, together with the dimension vector, determines the \'etale structure of the moduli space around the representation $V$.
\begin{thm}
For every point $V \in M^{ss}_\alpha(Q,\theta)$ we have an \'etale isomophism between an open neighborhood of the zero representation in $\iss_\beta Q_V$ and and open neighborhood of $V$.
\end{thm}

We now have almost everything we need to determine which moduli spaces are smooth projective varieties. We now only need to know the $\theta$-(semi)stable representations of our quiver. These were determined by Adriaenssens and Le~Bruyn in \cite{AdLeB}.
\begin{thm}
\begin{enumerate}
\item[]
\item
For a dimension vector $\alpha=(a_1,\dots,a_p;b_1,\dots,b_q)$ such that $\theta\cdot \alpha=0$ there exists always $\theta$-semistable representations, in this case we denote $n=\sum_{i=1}^p a_i$.
\item
$M^{ss}_\alpha(Q,\theta)$ contains a non-empty subset of $\theta$-stable representations, which is then
a dense open subset, if and only if
$$
\forall i\le p,j \le q: a_i+b_j \le n~(**)
$$
unless $\forall i,j: a_i=b_j$.
or $n=1$ in which case $M^{ss}_\alpha(Q,\theta)$ is just a point.
\end{enumerate}
Because we use the condition $(**)$ quite often in the next section we will call a dimension vector satisfying $(**)$ \emph{almost simple}.
\end{thm}

\section{Smoothness of $\Mss_\alpha(Q,\theta)$}
In this section we use the local quivers introduced earlier to determine which of the moduli spaces correspond to a smooth projective variety. 

Suppose that $\Mss_\alpha(Q,\theta)$ is smooth in the point corresponding to a representation
$$
V = S_1^{\oplus a_1}\oplus \dots \oplus S_k^{\oplus a_k}
$$
then the point $0$ must also be smooth in $\iss_\alpha Q_V$. For what kind of quivers $Q_V$ is this the case? 

If we look at the algebra of invariants $\C[\iss_\alpha Q]$, a well known theorem of Procesi and Le Bruyn \cite{ProLeB} states that this algebra is generated by a finite number of traces along cycles, $c_i$, modulo some relations.
$$
\C[\iss_\alpha Q] = \C[c_1,\dots,c_k]/(f_1,\dots f_l).
$$
This algebra inherits the grading of $\C[\Repr_\alpha Q]$ as the action of $GL(\alpha)$ preserves this grading. It is a well known fact that a positively graded, connected algebra is smooth if and only if it is a polynomial algebra (see for instance \cite{VON}).

Now we know the necessary condition for $\iss_\alpha Q$ one can try to classify all quivers and dimension vectors for which this $\iss_\alpha Q$  is indeed an affine space. Because this is a highly nontrivial problem we will only restrict ourselves to certain quivers with two vertices. These are the quivers that show up in the $\theta$-semistable representations that are a direct sum of two $\theta$-stables. Demanding that the moduli space is smooth in these of points will give us a restriction. The remaining cases we will consider in the next section and see that they are indeed totally smooth.

\begin{lemma}\label{lemma0}
The following quiver with indicated dimension vector has as ring of invariants a polynomial algebra if and only if there's at most one cycle connecting the two vertices (i.e. $k\le 1$)
$$
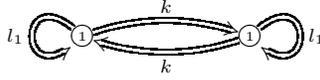
\xymatrix{
\vtx{1}\ar@/^/@2{->}[rr]^k\ar@(ul,dl)@2{->}[]_{l_1} &&\vtx{1}\ar@/^/@2{->}[ll]^k \ar@(ur,dr)@2{->}[]^{l_1}}
$$
\end{lemma}

\begin{prf}
The representation space is spanned by all loops $L_i$ in the both vertices and all cycles
$$
X_{ij} = a_ib_j
$$
All these cycles are neccesary to generate the algebra, because the representation for which all the arrows are zero except $a_i$ and $b_j$ is not equivalent to the zero and has as values in the cycles all zero's except for $X_{ij}$. The relations between the cycles are of the form
$$
X_{ij}X_{kl}=X_{il}X_{kj}
$$
These relations prevent $\iss_\alpha Q$ from being an affine space. The only way to make $\iss_\alpha Q$ into an affine space is that there's only $1$ such cycle.
\eop \end{prf}

If $M^{ss}_\alpha(Q,\theta)$ is a smooth space, it will be definitely smooth in the semisimple points which have only two factors with multiplicity $1$. We will see that this is not the case for most of the moduli spaces. By the previous lemma we only have to check that the number of arrows connecting the both factors is not greater than $1$, i.e.
$$
\chi_Q(\alpha_1,\alpha_2) \ge -1.
$$
using this fact we can deduce the following facts
\begin{lemma}\label{lemma1}
Suppose $\alpha = (a_1,\dots,a_p;b_1,\dots,b_q)$ is a simple dimension vector, then the degeneration 
$$(a_1,\dots,a_i-1,\dots,a_p;b_1,\dots,b_j-1,\dots,b_q)+\epsilon_{ij}$$
is smooth if and only if $a_i +b_j= n$. (In this degeneration $\epsilon_{ij}$ is short hand for the dimension vector $(\delta_{1i},\dots,\delta_{pi};\delta_{1j},\dots,\delta_{qj}$)
\end{lemma}
\begin{prf}
If we calculate the Euler form
$$
\chi(\alpha',\epsilon_{ij}) = a_i-1 + b_i-1 -(n-1) =-1 +(a_i+b_i-n)
$$
equals $-1$ if and only if $a_i+b_i=n$
\eop \end{prf}
In the folowing we suppose that the dimension vector is ordered i.e.
$$
a_1\ge \dots\ge a_p,~b_1\ge \dots\ge b_q
$$
\begin{lemma}\label{lemma2}
If $\alpha$ is an almost simple dimension vector and $a_1=a_2$ and $b_1=b_2$ and $a_1+b_1=n$ then $a_i=b_i=0,$ $i>2$
\end{lemma}
\begin{prf}
We know that $\sum a_i = n$ and $\sum b_i=n$ so
$$
\sum_{i=1}^p a_i +\sum_{j=1}^q b_j = a_1 + b_1 + a_2 +b_2 +\sum_{i>2} a_i +\sum_{j>2} b_j= 2n + \sum_{i>2} a_i +\sum_{j>2} b_j,
$$
Which implies that the last two terms must be zero and $a_1=a_2=b_1=b_2=\frac n2$.
\eop \end{prf}

\begin{lemma}\label{lemma3}
Suppose $\alpha = (a_1,\dots,a_p;b_1,\dots,b_q)$ is a dimension vector of a $\theta$-stable for which all the possible degenerations
$$
(a_1,\dots,a_i-1,\dots,a_p|b_1,\dots,b_j-1,\dots,b_q)+\epsilon_{ij}
$$
are smooth then either 
\begin{itemize}
\item $\alpha =(1;1)$, wich is the trivial case
\item $b_1=\dots=b_q$, $a_1>a_2$, and $a_1+b_1 =n$  or vice versa changing the $a's$ in the $b's$
\end{itemize}
\end{lemma}
\begin{prf}
Suppose we're not in the trivial case. If $a_1+b_1<n$  then we can choose whatever $\epsilon_{ij}$ we like and $\alpha-\epsilon_{ij}$ will be simple, but by the first lemma this degeneration will not be smooth. So $a_1+ b_1=n$.

If $a_1=a_2$ and $b_1=b_2$ then by the second lemma the dimension vector is of the form $(a,a;a,a)$ and doesn't correspond to a $\theta$-stable.

So suppose that $a_1>a_2$ then we will prove that all the $b_j$ will be equal. Indeed, if this were not the case then $b_i<b_1$. But than we can split of $\epsilon_{1i}$ to obtain a valid degeneration because $a_1+b_i<a_1+b_1=n$ this degeneration will not be smooth.
\eop \end{prf}

\begin{lemma}\label{lemma4}
Suppose $\alpha = (a_1,\dots,a_p;b_1,\dots,b_q)$, $p\le q$ is dimension vector of a $\theta$-stable for which all the possible degenerations in two different simple components are smooth then either 
\begin{itemize}
\item
$\alpha = (q-1,1|1,\dots,1)$
\item
$\alpha = (b,b|b-1,1,1)$
\item
$\alpha = (4,2|2,2,2)$
\end{itemize}
\end{lemma}
\begin{prf}
Suppose that $a_2=a_1-l,~l>0$. By Lemma \ref{lemma3} we know that all the $b_i$ must be equal. We now distinguish the folowing cases
\begin{itemize}
\item
If $l\le q-2$ then splitting of $\epsilon_{1j}$ for $1\leq j\leq l$ yields a term
$$
(a_1-l,a_1-l,\dots,a_p;\underbrace{b-1,\dots,b-1}_{l},b,\dots,b)
$$
which is an almost simple dimension vector which satisfies de conditions of Lemma \ref{lemma2} so $a_3 =0$ and $b=1$. this gives us $\alpha = (q-1,1;1,\dots,1)$ (possibility 4).
\item 
If $l=q-1$ then $b$ cannot be $1$ otherwise 
$$
a_1 + b=n~ \implies ~a_1 = q-1 ~\implies ~a_1-l=q-1 - (q-1) = 0
$$
which is impossible because $\theta(\alpha) = 0$. If $b\ge 2$ and $a_3$ is not zero then
$$
(a_1-l,a_1-l, a_3-1,\dots,a_p;b-1,\dots,b-1)
$$
is an almost simple dimension vector which satisfies the conditions of Lemma \ref{lemma2} so $q=2$ and we find the solution $(b,b-1,1;b,b)$. If $a_3=0$ then $2a_1-q+1=qb$ en $a_1=(q-1)b$ so $q=3$ and because
$$
(a_1-3,a_1-3;b-1,b-1,b-2)
$$
is an almost simple dimension vector which satisfies the conditions of lemma 2, therefore $b$ must be two and we obtain $(4,2;2,2,2)$.
\item
If $q\le l$ then we can split $\alpha$ in the folowing way:
\se{
\alpha = &(a_1-q+1,a_1-l-1,a_3,\dots,a_p;b_1-1,\dots,b_q-1) \\&+(q-1,1;b-1,\dots,b-1)
}
the Euler product for this degeneration is:
\se{
\chi &= (a_1-q+1)(q-1) + (a_1-l-1) +q(b-1)- q^2(b-1)\\
&=(q-1)^2(b-1) + ((q-1)b-l-1)+ q(b-1) - q^2(b-1)\\
&=(q-1)^2(b-1) + (2q-1)(b-1)+q-1-l-1 - q^2(b-1)\\
&=q-1-l-1<-1
}
which implies that it is not smooth.
\end{itemize}
\eop \end{prf}
We will now look at these quivers and determine which of them are really smooth in all their degenerations.
\begin{lemma}
The quiver with dimension vector $(b,b;b,b-1,1)$ is non-smooth for $b\geq 3$.
\end{lemma}
\begin{prf}
Consider the dimension vector for $b=3$, and look at the degeneration
$$(1,0;0,0,1)\oplus (1,0;0,1,0)^{\oplus 2}\oplus (0,1;1,0,0)^{\oplus 3}.$$
Its local quiver is
$$\xymatrix{\vtx{1}\ar@/^.3pc/[r] & \vtx{3}\ar@/^.3pc/[r]\ar@/^.3pc/[l] & \vtx{2}\ar@/^.3pc/[l]},$$
which has a non-smooth moduli space, as one of us has shown in \cite{Raf}.

Now consider the case when $b>3$ and look at the degeneration
$$(1,1;1,0,1)\oplus (1,1;1,1,0)^{\oplus (b-1)}.$$
Its local quiver becomes
$$
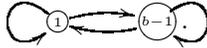
\xymatrix{\vtx{1}\ar@/^.3pc/[r]\ar@(ul,dl) & \vtx{b-1}\ar@/^.3pc/[l]\ar@(ur,dr)}.$$
This local quiver becomes non-smooth for the degeneration
$$
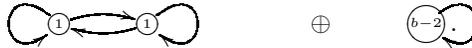
\xymatrix{\vtx{1}\ar@/^.3pc/[r]\ar@(ul,dl) & \vtx{1}\ar@/^.3pc/[l]\ar@(ur,dr) &&\oplus & \vtx{b-2}\ar@(ur,dr)}.$$
Indeed, the local quiver of this degeneration becomes
$$
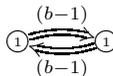
\xymatrix{\vtx{1}\ar@{=>}@/^.3pc/[r]^{(b-1)} & \vtx{1}\ar@{=>}@/^.3pc/[l]^{(b-1)}}$$
which is non-smooth according to Lemma \ref{lemma0}.
\eop \end{prf}

\section{Determining the structure of the moduli spaces}
In this section we will work out the structure of the moduli spaces associated to the quivers with dimension vectors as appearing in Lemma \ref{lemma4}.

First, we will consider is the quiver
$$
Q_m :=
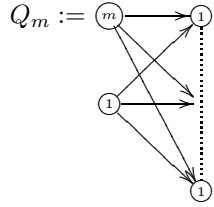
\xymatrix
{
\vtx{m} \ar[r]\ar[rd]\ar[rdd] & \vtx{1}\ar@{.}[dd]\\
\vtx{1} \ar[ru]\ar[r]\ar[rd]  & \\
\                             & \vtx{1}
}
$$
We will denote by $k_i$ (resp.~$c_i$) the arrow running from the first (resp.~second) vertex in the left part of the bipartite quiver to the $i$th arrow in the right part of the quiver
\begin{thm}
$\iss_\alpha(Q_m)$ is the projective space in $m$ dimensions.
\end{thm}
\begin{prf}
To prove the above statement we just have to prove that the ring of semi-invariants is the polynomial ring in $m+1$
variables. We first prove that this ring is generated by $m+1$ semi-invariants. 

All the semi-invariants are generated by the matrix-semi-invariants. Suppose now that we have a representation where the arrows $k_i$ are represented by row vectors $K_i$ and the arrows $c_i$ by constants $C_i$. A general matrix-semi-invariant of the order $l$ is of the form
$$
\left|\begin{matrix}
s_{11}K_1 & s_{12}C_1&\dots &s_{1,2l-1}K_1&s_{1,2l}C_1\\
\vdots&&&&\vdots\\
s_{m+1,1}K_{m+1} & s_{m+1,2}C_{m+1}&\dots &s_{m+1,2l-1}K_{m+1}&s_{m+1,2l}C_{m+1}\\
\vdots&&&&\vdots\\
s_{lm+l,1}K_{m+1} & s_{lm+l,2}C_{m+1}&\dots &s_{lm+l,2l-1}K_{lm+l}&s_{lm+l,2l}C_{m+1}\\
\end{matrix}\right|,
$$
where the $s_{ij}$ represent complex numbers. Using the multilinearity in the rows, one can rewrite the big determinant as a linear combination of determinants where on each row there's exactly one $s_{ij}$ equal to $1$ and all the others are zero.

Switching rows enables us to put all the rows where one of the two first $s's$ are non-zero above (take care to switch only rows modulo $m+1$). The number of such rows must be equal to $m+1$ otherwise the determinant will be zero.
$$
\xymatrix@R-1.5pc@C-1.5pc
{
&              &\ar@{.}[dddd]& &\\
&k\times m+1   &             & {\begin{matrix}0&\cdots& 0\end{matrix}}&\\
\ar@{.}[rrrr] &&&&\\
&{\begin{matrix} 0\\ \vdots\\ 0\end{matrix}} &             &{\begin{matrix}*&\cdots& *\\\vdots&&\vdots\\ *&\cdots& *\end{matrix}}&\\
&              &             & &
\save "2,1"."5,1"*+<-.1pc>\frm{(}
\restore
\save "2,5"."5,5"*+<-.1pc>\frm{)}
\restore
}
$$
So in the above matrix the upper left corner is a square $m+1 \times m+1$ dimensional matrix. The big deternimant now decomposes in a product of a semi-invariant of degree $1$ and one of degree $l-1$. By induction all the semi-invariants are generated by generated by the ones of degree $1$. When we take a look at these we can see that by the multilineary of the determinant every such semi-invariant is a linear combination of the folowing one's
$$
T_i := \left|
\begin{matrix}
K_1&0\\
\vdots&\vdots\\
K_{i-1}&0\\
0& C_i\\
K_{i+1}&0\\
\vdots&\vdots\\
K_{m+1}&0\\
\end{matrix}
\right|,~i=1 \to m+1
$$
Between those $m+1$ semi-invariants are no relations because if we consider an $m+1$-tuple different from zero say $(x_1,\dots,x_{m+1})$.
The representation
\se{
K_i &:= \begin{pmatrix}\delta_{1i}&\cdots&\delta_{m i}\end{pmatrix},~1\le i\le m\\
K_{m+1} &:= \begin{pmatrix}1&\cdots& 1\end{pmatrix}\\
C_{i}&:= x_i
}
has as invariant $T_i=x_i$. So we have a set of $m+1$ independant generators and hence the ring of semi-invariants is $\C[T_1t,\dots,T_{m+1}t]$.
\eop \end{prf}

We can apply a similar trick to 
$$
Q_{III} :=
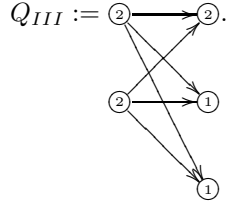
\xymatrix
{
\vtx{2}\ar[r]\ar[rd]\ar[rdd] & \vtx{2}\\
\vtx{2}\ar[ru]\ar[r]\ar[rd] & \vtx{1}\\
 & \vtx{1}
}.
$$
We denote the arrows between the vertices with dimension $2$ as $a$ and $b$ and the arrows from the first (resp.~second) vertex of the left part of the quiver to the vertices with dimension $1$ in the right part with $c_1$ and $c_2$ (resp.~$d_1$ and $d_2$).
\begin{thm}
$\iss_\alpha(Q_{III})$ is the projective space in $3$ dimensions.
\end{thm}
\begin{prf}
Suppose that we have a representation where the arrows $a$ and $b$ are represented by $2\times 2$-matrices $A,B$ and the arrows $c_i$, $d_i$ by row-vectors $C_i,D_i$. A general matrix-semi-invariant of the order $l$ is of the form
$$
\left|\begin{matrix}
s_{11}A & s_{12}B&\dots &s_{1,2l-1}A&s_{1,2l}B\\
s_{21}C_1 & s_{12}D_1&\dots &s_{2,2l-1}C_1&s_{2,2l}D_1\\
s_{31}C_2 & s_{12}D_2&\dots &s_{3,2l-1}C_2&s_{3,2l}D_2\\
\vdots&&&&\vdots\\
s_{3l, 1}C_2 & s_{12}D_2&\dots &s_{2,2l-1}C_2&s_{2,2l}D_2\\
\end{matrix}\right|,
$$
where the $s_{ij}$ represent complex numbers. Using the multilinearity in the rows, one can rewrite the big determinant as a linear combination of determinants where on each row of the $C,D$-part there's exactly one $s_{ij}$ equals $1$ and all the others are zero. For the $A,B$-part this isn't possible because they consist of two rows. But by subtracting and switching columns we can get a couple of rows of the form
$$
\begin{pmatrix}
s_1A &s_2B&0&\dots&0
\end{pmatrix}
$$
Switching rows enables us to put all the rows where one of the two first $s's$ are non-zero above (take care to switch only rows modulo $4$). The number of such rows must be equal to $4$ otherwise the determinant will be zero. As in the previous theorem the big determinant now decomposes in a product of a semi-invariant of degree $1$ and one of degree $l-1$. By induction all the semi-invariants are generated by generated by the ones of degree $1$. We take a look at these we can see that by the multilineary of the determinant every such semi-invariant is a linear combination of the folowing ones
$$
T_1 := \left|
\begin{matrix}
A&0\\
0& D_1\\
0&D_2
\end{matrix}\right|,
~T_2 := \left|
\begin{matrix}
0&B\\
C_1&0\\
C_2&0
\end{matrix}\right|,
~T_3 := \left|
\begin{matrix}
A&B\\
C_1&0\\
0&D_2
\end{matrix}\right|,
~T_4 := \left|
\begin{matrix}
A&B\\
0&D_1\\
C_2&0
\end{matrix}\right|.
$$

Between those $4$ semi-invariants are no relations because if we consider an $4$-tuple different from zero say $(x_1,\dots,x_{4})$.
The representation
\se{\begin{matrix}
A= \begin{pmatrix}
1&0\\
0&x_1
\end{pmatrix}&
B= \begin{pmatrix}
x_2&0\\
0&1
\end{pmatrix}
\\
C_1=\begin{pmatrix}
0&1
\end{pmatrix}&
C_2=\begin{pmatrix}
0&x_3
\end{pmatrix}\\
D_1=\begin{pmatrix}
1&0
\end{pmatrix}&
D_2=\begin{pmatrix}
-x_4&0
\end{pmatrix}
\end{matrix}}

has as invariant $T_i=x_i$. So we have a set of $m+1$ independent generators and hence the ring of semi-invariants is $\C[T_1t,\dots,T_{4}t]$.
\eop \end{prf}

We now have one more situation to look at.
\begin{thm}
The quiver 
$$
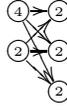
\xymatrix@R-1.5pc@C-1.5pc{
\vtx{4}\ar[r]\ar[dr]\ar[ddr] & \vtx{2}\\
\vtx{2}\ar[ur]\ar[r]\ar[dr]  & \vtx{2}\\
                             & \vtx{2}\\
}$$
with dimension vector $(4,2;2,2,2)$ has $\PP^5$ as its moduli space.
\end{thm}
\begin{prf}
Using the fact that the map 
$V_{\xi^1}\rightarrow W_{\eta^1}\oplus W_{\eta^2}\oplus W_{\eta^3}$
must be injective, we can apply \emph{reflection functors} to identify the moduli space of the original quiver with the moduli space of the reflected quiver
$$
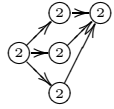
\xymatrix@R-1.5pc@C-1.5pc{
                             & \vtx{2}\ar[r]& \vtx{2}\\
\vtx{2}\ar[ur]\ar[r]\ar[dr]  & \vtx{2}\ar[ur]&\\
                             & \vtx{2}\ar[uur]&\\
}$$
Which on its turn, using facts from invariant theory may be identified with the moduli space of the quiver
$$\xymatrix{
\vtx{2}\ar[r]\ar@/^/[r]\ar@/_/[r] & \vtx{2}
}$$
for $\theta = (-1,1)$. By results of Barth \cite{Ba} and Hulek \cite{Hu}, this moduli space indeed is a $\PP^5$.
\eop \end{prf}

Summarizing all results obtained in this paper we achieve the following theorem.
\begin{thm}
For the quiver $Q$ the only dimension vectors for which $\Mss_\alpha(Q,\theta)$ is smooth are in fact
\begin{itemize}
\item $\alpha := (m,1;1,\dots, 1)$ for which $\Mss_\alpha(Q,\theta)=\PP^m$.
\item $\alpha := (2,2;2,1, 1)$ for which $\Mss_\alpha(Q,\theta)=\PP^3$.
\item $\alpha := (4,2;2,2,2)$ for which $\Mss_\alpha(Q,\theta)=\PP^5$.
\end{itemize}
\end{thm}

\bibliography{biparbib}
\end{document}